\documentclass[12pt,reqno]{amsart}

\usepackage[latin1]{inputenc}
\usepackage{graphics,color}
\usepackage{amssymb,amsmath,amsthm,amscd}
\usepackage{latexsym,verbatim,graphicx,amsfonts}
\usepackage{hyperref}
\usepackage[mathscr]{euscript}
\usepackage{amsmath, amsthm, amssymb}
\usepackage{dsfont}
\usepackage{enumitem}
\usepackage{mathrsfs}
\usepackage{float}

\theoremstyle{plain}
\newtheorem{theorem}{Theorem}[section]

\newtheorem{corollary}[theorem]{Corollary}
\newtheorem{lemma}[theorem]{Lemma}
\newtheorem{proposition}[theorem]{Proposition}
\theoremstyle{definition}

\newtheorem{definition}[theorem]{Definition}

\theoremstyle{remark}
\newtheorem{rem}[theorem]{Remark}
\newtheorem{rems}[theorem]{Remark}

\newtheorem{question}[theorem]{Question}

\DeclareMathOperator{\spn}{span}

\newcommand{\norm}[1]{\left\lVert#1\right\rVert}

\DeclareMathOperator{\supp}{supp}
\DeclareMathOperator{\identity}{Id}
\DeclareMathOperator{\real}{Re}

\DeclareMathOperator{\proj}{Proj}

\newcommand{\finset}[1]{\left[{#1}\right]^{< \omega}}
\newcommand{\inset}[1]{\left[{#1}\right]^{\leq \omega}}
\newcommand{\id}[1]{\mathds{1}_{#1}}

\newcommand{\vertiii}[1]{{\left\vert\kern-0.25ex\left\vert\kern-0.25ex\left\vert #1 
    \right\vert\kern-0.25ex\right\vert\kern-0.25ex\right\vert}}

\newcommand{\Tsupp}[1]{\supp \left(T^*{#1} \right)}
\newcommand{\nonvan}[1]{\mathbf{Z}^c\left({#1}\right)}

\font\sstext=ecss1000
\font\sssub=ecss1000 at 7pt
\font\sssubsub=ecss1000 at 5pt

\newfam\ssfam
\textfont\ssfam=\sstext
\scriptfont\ssfam=\sssub
\scriptscriptfont\ssfam=\sssubsub

\def\ss{\fam\ssfam\sstext} 
\def\CH{{\ss CH}} 

\subjclass{46H10,
47L10
(primary); 
46E15,
47B01,
47L20 (secondary).}
\keywords{Banach space, bounded operator, closed operator ideal, quotient algebra, uniqueness of algebra norm, quantitative factorization of operators, $C(K)$-spaces, Calkin algebra, compact operators.}

\begin{document}

\title[Factorizations and minimality]{Factorizations and minimality of the Calkin Algebra norm for $C(K)$-spaces}
\author[A. Acuaviva]{Antonio Acuaviva}
\address{School of Mathematical Sciences,
Fylde College,
Lancaster University,
LA1 4YF,
United Kingdom} \email{ahacua@gmail.com}

\date{\today}

\begin{abstract}
For a scattered, locally compact Hausdorff space $K$, we prove that the essential norm on the Calkin algebra \break $\mathscr{B}(C_0(K))/\mathscr{K}(C_0(K))$ is a minimal algebra norm. The proof relies on establishing a quantitative factorization for the identity operator on $c_0$ through non-compact operators $T: C_0(K) \to X$, where $X$ is any Banach space that does not contain a copy of $\ell_1$ or whose dual unit ball is weak$^*$ sequentially compact. It follows that, for every ordinal $\alpha$, the algebras $\mathscr{B}(C[0,\alpha]))$ and $\mathscr{B}(C[0,\alpha]))/\mathscr{K}(C[0,\alpha]))$ have an unique algebra norm.
\end{abstract}

\maketitle

\section{Introduction and Main Results}\label{Introduction}

\subsection{Algebra norms, minimality and maximality of norms.}

By an algebra norm on a (real or complex) algebra $\mathcal{A}$ we mean a vector norm on $\mathcal{A}$ which is submultiplicative. We say that an algebra $\mathcal{A}$ admits an unique norm if any two algebra norms on $\mathcal{A}$ are equivalent. We also have the weaker notion of admitting a unique \emph{complete} norm, which simply says that any two complete algebra norms on $\mathcal{A}$ are equivalent. \\

Central to uniqueness-of-norm problems are the notions of minimality and maximality of a norm. An algebra norm $\norm{\cdot}$ on $\mathcal{A}$ is minimal (respectively, maximal) if for any other algebra norm $\vertiii{\cdot}$ on $\mathcal{A}$, there exists a constant $c > 0$ such that $c \norm{\cdot} \leq \vertiii{\cdot}$ (respectively, there exists $C > 0$ such that $\vertiii{\cdot} \leq C \norm{\cdot}$). Observe that, if an algebra norm is both minimal and maximal, then it is the unique norm for the algebra $\mathcal{A}$. Similarly, if a complete norm is either minimal or maximal, then the open mapping theorem gives that it is the unique complete norm. 

Therefore, uniqueness-of-norm problems can be reframed as problems about the minimality and maximality of norms. These two problems are, in general, quite different. We shall focus on these problems for the Banach algebra $\mathscr{B}(X)$ of bounded operators on a Banach space $X$, as well as its quotients. We start by examining the problem of minimality of the norm. \\

Among the first results in this direction is a celebrated theorem of Eidelheit which states that, for any Banach space $X$, the operator norm on the algebra $\mathscr{B}(X)$ is minimal \cite{eidelheit1940isomorphisms} (and therefore it is the unique complete norm). A natural next step is to ask whether similar phenomena occur for quotients of $\mathscr{B}(X)$, most notably the Calkin algebra $\mathscr{B}(X)/\mathscr{K}(X)$, where $\mathscr{K}(X)$ denotes the ideal of compact operators. In this direction, Meyer \cite{meyer1992topological} showed minimality of the essential norm for $X = c_0$ and $X = \ell_p$, $1 \leq p < \infty$. In his PhD dissertation, Ware undertook a systematic study of the uniqueness-of-norm problem for Calkin algebras, significantly extending Meyer's work. He proved the minimality of the essential norm for a large class of Banach spaces $X$, including finite direct sums of $c_0$ and $\ell_p$ spaces, the infinite sums $\left( \bigoplus_{n \in \mathbb{N}} \ell_p^n \right)_{c_0}$ and $\left( \bigoplus_{n \in \mathbb{N}} \ell_p^n \right)_{\ell_q}$, Tsirelson's space $T$ and its dual $T^*$, and the James spaces $J_p$ for $1 < p < \infty$ \cite{ware2014uniqueness}. More recently, Johnson and Phillips \cite{johnsonandphillips} have shown the minimality of the essential norm in the Calkin algebra for $X = L_p$, $1 < p < \infty$, while Laustsen and Arnott \cite{arnott2023uniqueness} have shown minimality whenever $X = C_0(K_\mathcal{A})$, where $K_{\mathcal{A}}$ is a Mr\'owka space. For a broader and more in-depth perspective on this problem, we refer the reader to the recent paper by Arnott and Laustsen \cite{arnott2023uniqueness}. For a more detailed discussion specifically in the context of Calkin algebras, we recommend the PhD thesis of Ware \cite{ware2014uniqueness} and the survey article by Skillicorn \cite{skillicorn2015uniqueness}. 

Our main result for the minimality problem is summarised in the following theorem.

\begin{theorem}\label{th: main-minimality}
    Let $K$ be a scattered, locally compact Hausdorff space and consider the Calkin algebra $\mathscr{B}(C_0(K))/\mathscr{K}(C_0(K))$. Then the essential norm $\norm{\cdot}_e$ is minimal. In particular, the Calkin algebra admits a unique complete algebra norm.
\end{theorem}

The proof of the previous result is based on a bounding technique inspired by the Eidelheit-Yood method, as elegantly presented in Ware's thesis \cite[Theorem 1.1.2] {ware2014uniqueness}, together with a result about uniform factorizations of the identity on $c_0$, Theorem \ref{th: fact-iden}, which we present later and believe may be of independent interest. \\

We now turn to the question of the maximality of norms. To highlight the contrast with the minimality problem, observe that while Eidelheit \cite{eidelheit1940isomorphisms} proved that the operator norm is minimal in $\mathscr{B}(X)$ for any Banach space $X$, an example of Read \cite{read1989discontinuous} gives a Banach space such that the operator norm on $\mathscr{B}(X)$ is not maximal (equivalently, there are spaces $X$ such that $\mathscr{B}(X)$ admits a discontinuous algebra homomorphism into a Banach algebra). \\

Observe that the maximality of the norm passes down to quotient spaces, while minimality does not. This is because the maximality of the norm is equivalent to proving that every homomorphism from the algebra into a Banach algebra is continuous. This characterisation explains why, when proving the maximality of the norm on a Calkin algebra, one usually focuses on proving the same question for the algebra of operators $\mathscr{B}(X)$. This gives rise to the following natural question.

\begin{question}\label{qs: 1}
     Is there a scattered, locally compact Hausdorff space $K$ such that $\mathscr{B}(C_0(K))$ admits a discontinuous homomorphism into a Banach algebra?
\end{question}

Under the continuum hypothesis (\CH), our previous question admits a positive answer if one drops the scattered condition, see \cite[Remark 40]{koszmider2021banach}, however, the space discussed there is far from being scattered. In upcoming work, we have managed to construct a scattered, locally compact Hausdorff space $K$ such that, under (\CH), $\mathscr{B}(C_0(K))$ admits a discontinuous homomorphism into a Banach algebra \cite{Acuaviva1}. The question remains open if the assumption on (\CH) is dropped. \\

When dealing with this automatic continuity question for operator algebras $\mathscr{B}(X)$, a famous theorem of B. E. Johnson \cite{johnson1967continuity} shows that it is enough for the space $X$ to be isomorphic to its square $X \oplus X$ (or, more generally, to admit a \textit{continuous bisection of the identity}).
Nevertheless, this criterion alone proves insufficient for numerous spaces, even within the realm of continuous functions on scattered compact Hausdorff spaces. As an illustration, Semadeni \cite{semadeni1960banach} showed that $C([0, \omega_1])$ fails to be isomorphic to its square, where $\omega_1$ represents the first uncountable ordinal. Even more, Loy and Willis \cite{loy1989continuity} went on to prove that this space cannot admit a continuous bisection of the identity, thus rendering Johnson's result inapplicable.  \\

That being said, Ogden \cite{ogden1996homomorphisms} has proved that every homomorphism from $\mathscr{B}(C([0, \omega_\eta]))$, where $\omega_\eta$ is the first ordinal of cardinality $\aleph_\eta$ for an ordinal $\eta$, to a Banach algebra is continuous. We provide a brief account of how Ogden's result extends to all ordinals $\alpha$, yielding the following.

\begin{theorem}\label{th: Ogden}
    Let $\alpha$ be an ordinal. Then any homomorphism from $\mathscr{B}(C([0, \alpha]))$ into a Banach algebra is continuous.
\end{theorem}

The previous theorem, coupled with Theorem \ref{th: main-minimality} and Eidelheit's theorem, gives a solution to the uniqueness of norm problem for the algebra of operators and the Calkin algebra on ordinal spaces.

\begin{corollary}
    Let $\alpha$ be an ordinal. Then the algebra of operators $\mathscr{B}(C([0, \alpha]))$ and the Calkin algebra $\mathscr{B}(C([0, \alpha]))/\mathscr{K}(C([0, \alpha]))$ admit an unique algebra norm.
\end{corollary}

Lastly, we turn our attention to another property of interest in the study of algebra norms, namely \emph{Bonsall's minimality property}, see \cite{Tylli1}. The formal definition is as follows, an algebra norm $\norm{\cdot}$ on an algebra $\mathcal{A}$ has Bonsall's minimality property if $\norm{\cdot} = \vertiii{\cdot}$ whenever $\vertiii{\cdot}$ is any other algebra norm on $\mathcal{A}$ satisfying $\vertiii{\cdot} \leq \norm{\cdot}$. We obtain the following.

\begin{theorem}\label{th: main-bonsall-min}
    Let $K$ be a scattered, locally compact Hausdorff space and consider the Calkin algebra $\mathscr{B}(C_0(K))/\mathscr{K}(C_0(K))$. Then the dual essential norm $T \mapsto \norm{T^*}_e$ has Bonsall's minimality property.
\end{theorem}

\subsection{Quantitative factorizations of the identity.} 

Factorizations of idempotent operators play a prominent role in the study of the algebra $\mathscr{B}(X)$ of bounded operators on a Banach space $X$. For instance, such factorizations are central to the classification of the closed ideals of $\mathscr{B}(X)$, and they also feature prominently in the study of the uniqueness of algebra norms. In our case, we will establish the following quantitative result.

\begin{theorem}\label{th: fact-iden}
    Let $K$ be a scattered, locally compact Hausdorff space and let $X$ be a Banach space that does not contain a copy of $\ell^1$ or for which the unit ball of $X^*$ is weak$^*$ sequentially compact. Then for any non-compact operator $T: C_0(K) \to X$ and any $\varepsilon > 0$, there exist operators $U \in \mathscr{B}(X, c_0)$ and $V \in \mathscr{B}(c_0, C_0(K))$ such that
    \begin{equation*}
        UTV = I_{c_0} \hspace{5pt} \text{ and } \hspace{5pt} \norm{U} \norm{V} < \frac{2}{\norm{T}_e} + \varepsilon.
    \end{equation*}
\end{theorem}

We also provide a further improvement of the previous theorem, namely, we show how the upper bound on the norms of the operators $U$ and $V$ is intimately related to the essential norm of the adjoint operator of $T$, obtaining a characterisation of this norm.

\begin{theorem}\label{th: ess-norm-adjoint} 
     Let $K$ be a scattered, locally compact Hausdorff space and let $X$ be a Banach space that does not contain a copy of $\ell^1$ or for which the unit ball of $X^*$ is weak$^*$ sequentially compact. Then for any non-compact operator $T: C_0(K) \to X$ 
     \begin{equation*}
         \norm{T^*}_e = \sup \{\norm{U}^{-1}\norm{V}^{-1}: UTV = I_{c_0} \}.
     \end{equation*}
\end{theorem}

\begin{rem}
    Theorem \ref{th: fact-iden} could also be derived from Theorem \ref{th: ess-norm-adjoint} together with a duality result due to S. Axler, N. Jewell, and A. Shields \cite{axler1980essential}. Although their result is stated only for the case where the domain and codomain coincide, a careful reading of the proof shows that it suffices for the dual of the domain to have the $\lambda$-metric approximation property for their argument to go through.
\end{rem}
\begin{rem}
    If $K$ has an accumulation point, i.e. if $C_0(K) \not = c_0(\Gamma)$, the constant $2$ is sharp in Theorem \ref{th: fact-iden}, since in this case there exists an operator $T \in \mathscr{B}(C_0(K))$ such that $\norm{T}_e = 2$ while $\norm{T^*}_e = 1$. Clearly, if $C_0(K) = c_0(\Gamma)$ we have $\norm{T}_e = \norm{T^*}_e$.
\end{rem}

While our motivation regarding this factorization stems from the aforementioned paper of Arnott and Laustsen \cite{arnott2023uniqueness}, the origin of the ideas predates significantly in time. The starting point is the following classical theorem of Pe{\l}czy{\'n}ski. Observe that when $K$ is scattered, compactness and weak compactness coincide.

\begin{theorem}[Pe{\l}czy{\'n}ski \cite{pelczynski1962banach}]\label{th: pelc}
    Let $K$ be a locally compact Hausdorff space and $X$ be a Banach space. Then an operator $T: C_0(K) \to X$ is non-weakly compact if and only if $T$ fixes a copy of $c_0$.
\end{theorem}

We include a brief proof, based on the presentation in \cite[Theorem 5.5.3]{albiac2006topics}, as the reader may observe that the arguments in our proof of Theorem \ref{th: fact-iden} constitute quantitative refinements of those presented therein.

\begin{proof}

Let $T: C(K) \to X$ be non-weakly compact. By Gantmacher's theorem, $T^*: X^* \to C(K)^*$, where $C(K)^* = \mathcal{M}(K)$ is the space of all finite regular Borel measures on $K$, is also non-weakly compact; that is, the bounded set $T^*(B_{X^*}) \subseteq \mathcal{M}(K)$ is non-relatively weakly compact. By Grothendieck's classical characterisation of non-weak compactness in these spaces \cite{grothendieck1953applications}, there exist $\delta > 0$, a disjoint sequence of open sets $(U_n)_{n=1}^\infty$ in $K$, and a sequence $(x^*_n)_{n=1}^\infty \subseteq B_{X^*}$ such that if $\nu_n = T^*x^*_n$ then $\nu_n(U_n) > \delta$ for all $n \in \mathbb{N}$.

For each $n$ we can find a compact set $F_n \subseteq U_n$ such that $\nu_n(U_n \backslash F_n) < \delta/2$. By Urysohn's lemma, there exists $f_n \in C(K)$, $0 \leq f_n \leq 1$ such that $f_n = 0$ in $K \backslash U_n$ and $f_n = 1$ in $F_n$. The map $S: c_0 \to C(K)$ defined by $Se_n = f_n$ is an isometric embedding of $c_0$ into $C(K)$, where $(e_n)_{n=1}^\infty$ is the canonical basis of $c_0$.

Note that the operator $TS: c_0 \to X$ satisfies
\begin{equation*}
    \inf_{n \in \mathbb{N}} \norm{TSe_n} \geq \inf_{n \in \mathbb{N}}|x^*_n(TSe_n)| = \inf_{n \in \mathbb{N}} |T^*x^*_n(f_n)| \geq \delta/2.
\end{equation*}
It follows from a result of Rosenthal (see \cite{rosenthal1970relatively}, Remark 1 after Theorem 3.4), that $TS$ is an isomorphism on a copy of $c_0$, and thus the same holds for $T$.
\end{proof}

One can see from the proof that the value
\begin{equation*}
    \varepsilon(T) = \sup \left\{ \inf_n \left\| (T^* x_n^*)|_{U_n} \right\| : (U_n)_{n=1}^\infty \subseteq K, (x_n^*)_{n=1}^\infty \subseteq B_{X^*} \right\},
\end{equation*}
where the supremum is taken over all sequences of disjoint open sets $(U_n)_{n=1}^\infty \subseteq K$ and all sequences $(x_n^*)_{n=1}^\infty \subseteq B_{X^*}$, quantifies how well $T$ preserves a copy of $c_0$. For the factorization, observe that the stability conditions on the codomain, i.e. that $X$ contains no copy of $\ell_1$ or that the unit ball of $X^*$ is weak$^*$ sequentially compact, are there to guarantee that for any copy of $c_0$ we can pass down to a sub-copy that is complemented with good constants. \\

The essence of the proof of Theorem \ref{th: fact-iden} lies in finding a lower bound for $\varepsilon(T)$ in terms of the essential norm $\norm{T}_e$. Therefore, Theorem \ref{th: fact-iden} can be interpreted as a quantitative version of Theorem \ref{th: pelc}. Since the groundwork for deriving the quantitative factorization from this was already explored and developed by Arnott and Laustsen \cite[Section 7]{arnott2023uniqueness}, we focus on showing an equivalent condition, in the form of Theorem \ref{th: main} (see Section 3).

\begin{rem}
    In light of Pe{\l}czy{\'n}ski's theorem, one might hope that the assumption that $K$ is scattered in Theorem \ref{th: fact-iden} could be removed by replacing non-compactness with non-weak compactness. However, this is not possible, as $\varepsilon( \cdot )$ is not comparable in general with the weak essential norm $\norm{\cdot}_w$. Details of this will appear elsewhere.
\end{rem}

\bigskip
\section{Organization}\label{sec-organization}

We now outline the structure of the paper. Sections \ref{sec-sketch} to \ref{sec-sharpness and improvements} present the proofs concerning the quantitative factorization of the identity on $c_0$, while Sections \ref{sec-minimality} and \ref{sec-uniqueness-norm} address the minimality and maximality of algebra norms. \\

In Section \ref{sec-sketch} we introduce some background material and terminology. Following the work of Arnott and Laustsen, we provide a reformulation of Theorem \ref{th: fact-iden} in the form of Theorem \ref{th: main},  which we will use as our starting point for the proof. We provide a proof of Theorem \ref{th: fact-iden} subject to the verification of Theorem \ref{th: main}.

Sections \ref{sec-top-prel} to \ref{sec-main-proof} work through the proof of Theorem \ref{th: main}. In section \ref{sec-top-prel} we present some topological results, which allow us to introduce a particular class of projection operators. In the last part of this section, we will recall some properties of $C_0(K)^*$, which will be needed moving forward.  \\

Section \ref{sec-technical} establishes the two main ingredients for the proof. First, in Proposition \ref{prop-almost-orthonormal-system} we prove that we can find a sequence of functionals $(x^*_j)_{j=1}^\infty \subset B_{X^*}$ which witness the action of $T$ and that are almost orthonormal, in the sense that their weight is supported in a family of points $(F_j)_{j=1}^\infty$ which are disjoint.

The other ingredient deals with the accumulation of weight of an infinite collection of functionals. In particular, we establish that this weight cannot significantly accumulate on too many families of disjoint points. All of this is made precise in Lemma \ref{lmm-ess-weight}. Technical considerations make us restate this result in the form of Corollary \ref{corollary-ess-weight}. \\

Section \ref{sec-main-proof} contains the proof of Theorem \ref{th: main}. The key step is an application of Corollary \ref{corollary-ess-weight} and a diagonal argument to the family of functionals and points given by Proposition \ref{prop-almost-orthonormal-system}, which gives us a refinement of the functionals and points, with the added benefit that we can now separate the points via disjoint neighbourhoods. This is achieved in Lemma \ref{lmm: auxilary-main}. The proof of the theorem then follows naturally.

Section \ref{sec-sharpness and improvements} deals with the relation between the factorizations of the identity on $c_0$ and the essential norm of the adjoint operator. This section contains the proof of Theorem \ref{th: ess-norm-adjoint}. \\

Lastly, Section \ref{sec-minimality} focuses on the question of the minimality of the Calkin algebra norm, including a proof of Theorem \ref{th: main-minimality}, which is achieved via a simultaneous factorization of operators in the form of Lemma \ref{lmm: multiple-factori}. This section also includes a proof of Theorem \ref{th: main-bonsall-min}. In Section \ref{sec-uniqueness-norm}, we explain how to go from Ogden's result to Theorem \ref{th: Ogden}.
\bigskip
\section{Preliminaries and proof of Theorem \ref{th: fact-iden}}\label{sec-sketch}

All normed spaces and algebras are over the scalar field $\mathbb{K}$, either the real or complex numbers, and we adhere to standard notational conventions. For a locally compact Hausdorff space $K$, $C_0(K)$ denotes the Banach space of continuous functions $f: K \to \mathbb{K}$ which vanish at infinity, in the sense that $\{k \in K: |f(k)| \leq \varepsilon \}$ is compact for all $\varepsilon > 0$, and equip it with the supremum norm.

We write $B_X$ for the closed unit ball of a normed space $X$ and denote by $I_X$ (or $I$ when $X$ is clear from the context) the identity operator on $X$. 

The term ``operator" will refer to a bounded linear map between normed spaces. As usual, given two normed spaces $X$ and $Y$, we denote by $\mathscr{B}(X,Y)$ the collection of all operators from $X$ to $Y$; if $X = Y$, we simply denote it by $\mathscr{B}(X)$. Similarly, we denote by $\mathscr{K}(X, Y)$ the ideal of compact operators, and for brevity we write $\norm{T}_e$ for the essential norm of an operator $T \in \mathscr{B}(X,Y)$, namely
\begin{equation*}
    \norm{T}_e = \norm{T + \mathscr{K}(X,Y)} = \inf \{\norm{T - S }: S \in \mathscr{K}(X,Y)\}.
\end{equation*}
More specialized notation will be introduced as and when needed.

In their work, Arnott and Laustsen established the following quantitative factorization result for the identity on $c_0$ \cite[Theorem 3.1]{arnott2023uniqueness}, which we restate below with notation adapted to our context.

\begin{theorem}[Arnott and Laustsen]\label{thr-niels}
    Let $T \in \mathscr{B}(C_0(K), X)$ where $K$ is a locally compact Hausdorff space and $X$ is a Banach space for which the unit ball of $X^*$ is weak$^*$ sequentially compact, and let $\delta > 0$. 

    Then $C_0(K)$ contains a sequence $(\eta_n)_{n \in \mathbb{N}}$ such that
    \begin{equation}\label{eq: cond1}
        \sup_{k \in K} \sum_{n=1}^\infty |\eta_n(k)| \leq 1 \hspace{5pt} \text{ and } \hspace{5pt} \inf_{n \in \mathbb{N}} \norm{T\eta_n} > \delta
    \end{equation}
    if and only if there are operators $U \in \mathscr{B}(X, c_0)$ and $V \in \mathscr{B}(c_0, C_0(K))$ such that
    \begin{equation*}
        UTV = I_{c_0} \hspace{5pt} \text{ and } \hspace{5pt} \norm{U} \norm{V} < 1/\delta.
    \end{equation*}
\end{theorem}

In the previous theorem, we emphasise that $\delta$ serves as the inverse of an upper bound for $\norm{U} \norm{V}$. Consequently, an upper bound on $\delta$ provides control over how small the upper bound for $\norm{U} \norm{V}$ can be. 

In their proof, Arnott and Laustsen note that the hypotheses can be relaxed to also include the case where $X$ does not contain a copy of $\ell^1$, provided that a result of Galego and Plichko \cite[Theorem 4.3]{galego2003banach} can also be applied in the case of complex scalars. The only obstruction in the proof of Galego and Plichko is the verification of the following result of Hagler and Johnson \cite[Theorem 1.a]{hagler1977banach} in the case of complex scalars.

\begin{theorem}[Hagler and Johnson]\label{th: hagler}
    Let $X$ be a Banach space. If $X^*$ contains a closed subspace in which no normalized sequence converges weak$^*$ to zero, then $\ell^1$ is isomorphic to a subspace of $X$.
\end{theorem}

We thank Professor W.~B.~Johnson for kindly explaining how the complex case follows from the real version of the theorem, and for permitting us to include his argument here.

\begin{proof}
    Let $X_\mathbb{R}$ be $X$ when considered as a real Banach space. It is elementary that a sequence is weakly convergent (respectively, weakly Cauchy) in $X$ if and only if it is weakly convergent (respectively, weakly Cauchy) in $X_\mathbb{R}$. Therefore, by the complex version of Rosenthal's $\ell_1$ theorem (due to Dor), it is enough to show that $X_\mathbb{R}$ contains a copy of (real) $\ell_1$. The map $x^* \mapsto \real x^*$ is a real-linear, isometric, weak$^*$-to-weak$^*$ homeomorphism of $X^*$ onto $(X_\mathbb{R})^*$. Consequently, if the complex space $X$ satisfies the hypothesis of Theorem \ref{th: hagler}, then so does the real space $X_\mathbb{R}$. The real version of the theorem then gives that $X_\mathbb{R}$ contains a copy of (real) $\ell_1$ and thus, as previously explained, $X$ contains a copy of (complex) $\ell_1$.
\end{proof}

Taking all of this into account, Theorem \ref{th: fact-iden} will follow from the next theorem, whose proof we shall present in Sections \ref{sec-top-prel} to \ref{sec-main-proof}.

\begin{theorem}\label{th: main}
    Let $K$ be a scattered, locally compact Hausdorff space, $X$ be a Banach space, $T: C_0(K) \to X$ be an operator with $\norm{T}_e = 1$ and $\delta \in (0,1/2)$. Then there exists a sequence of disjoint functions $(\eta_j)_{j=1}^\infty \subset B_{C_0(K)}$ such that $\norm{T\eta_j} > \delta$ for all $j \in \mathbb{N}$.
\end{theorem}
\bigskip
\section{Topological preliminaries and $C_0(K)$}\label{sec-top-prel}

In this section, and throughout the rest of the paper, $K$ will denote a scattered, locally compact Hausdorff space. We observe that $K$ is zero-dimensional, that is, there exists a base of clopen sets. To see this, note that $K$ is scattered, hence totally disconnected, while from general topology we have that a locally compact Hausdorff space is zero-dimensional if and only if it is totally disconnected \cite[Proposition 3.1.7]{tkachenko2008topological}, thus our space $K$ is zero-dimensional. Furthermore, since $K$ is also locally compact, for every $k \in K$, we can find a local base formed by clopen compact sets.

From now on, whenever we refer to a neighbourhood, we will always do so with the implicit understanding that it is clopen and compact. We denote by $\mathcal{T}$ the collection of all compact and clopen subsets of $K$, which is a base for the topology of $K$. 

Now we move to the idea of separating a finite family of points by a family of disjoint neighbourhoods. We make this precise in the next definitions.

\begin{definition}
    A family of sets $\mathcal{A}$ is called \emph{pairwise disjoint} if $A \cap B = \emptyset$ for any distinct $A, B \in \mathcal{A}$.
\end{definition}

\begin{definition}\label{def: nhg-collection}
    Given a subset $F \subset K$, we say that $\mathcal{U} \subset \mathcal{T}$ is a \emph{neighbourhood collection around $F$} if there exists a bijection $U: F \to \mathcal{U}$ such that $k \in U(k)$ for each $k \in F$.
\end{definition}

Combining the previous definitions, we will say that $\mathcal{U}$ is a \emph{pairwise disjoint neighbourhood collection around $F$} if $\mathcal{U}$ is a neighbourhood collection around $F$ and it is pairwise disjoint. 

\begin{definition}
    We define the \emph{support} of $\mathcal{U} \subset \mathcal{T}$ as
    \begin{equation*}
         S(\mathcal{U}) = \bigcup_{U \in \mathcal{U}} U,
    \end{equation*}
    in other words, the set of those points $k \in K$ which belong to some member of $\mathcal{U}$. 
\end{definition}

Since we shall mostly work with finite or countable collections of points and neighbourhoods, the following notation will be useful.

\begin{definition}
    Let $A$ be any set, we denote $\finset{A} = \{B \subset A: |B| < \omega\}$ and $\inset{A} = \{B \subset A: |B| \leq \omega \}$, in other words, the collections of finite and countable subsets of $A$, respectively.
\end{definition}

Our next result is automatic given that the space $K$ is Hausdorff.

\begin{lemma}\label{lmm-neigh-collection}
    Let $F \in \finset{K}$, then there exists a pairwise disjoint neighbourhood collection around $F$.
\end{lemma}

Having introduced the necessary topological notions, we are ready to discuss an important class of operators on $C_0(K)$.

\begin{definition}
    Let $F \in \finset{K}$ and let $\mathcal{U} \in \finset{\mathcal{T}}$ be a pairwise disjoint neighbourhood collection around $F$. We define
    \begin{equation*}
        \proj(\mathcal{U}, F) = \sum_{k \in F} \delta_{k}  \otimes \id{U(k)},
    \end{equation*}
    where $\delta_k: C_0(K) \to \mathbb{K}$ is the evaluation functional at $k$, $\id{U}$ denotes the indicator function of $U$ and $\delta_{k}  \otimes \id{U(k)}$ denotes the rank-one operator defined via $f \mapsto f(k)\id{U(k)}$.
\end{definition}

Note that $\id{U} \in C_0(K)$ because it vanishes at infinity since $U$ is compact and it is continuous since $U$ is clopen. Hence $\proj(\mathcal{U}, F)$ defines an operator $C_0(K) \to C_0(K)$ with the following properties, whose verification is straightforward.

\begin{lemma}\label{lmm: proj}
        \begin{enumerate}[label={\normalfont{(\alph*)}}]
            \item \label{lmm: proj-a} $\proj(\mathcal{U}, F)$ is a finite-rank projection.
            \item \label{lmm: proj-d} $\norm{\proj(\mathcal{U}, F)} = 1$, except in the trivial case where $F = \emptyset$. Hence $\norm{I - \proj(\mathcal{U}, F)} \leq 2$.
            \item \label{lmm: proj-b} For any $k \in F$, the adjoint projection satisfies $\proj(\mathcal{U}, F)^* \delta_{k} = \delta_{k}$.
        \end{enumerate}
\end{lemma}

Finally, we shall need some properties of $C_0(K)^*$. We start by recalling a well-known variant of Rudin's famous theorem \cite{rudin1957continuous}. For completeness, we include a brief proof.

\begin{theorem}\label{thr-ell1}
    Let $K$ be a scattered, locally compact Hausdorff space. Then $C_0(K)^* \cong \ell^1(K)$. 
\end{theorem}

\begin{proof}
 Let $\alpha K$ be the one-point compactification of $K$, adding the point at $\infty$. Since $\alpha K$ is also scattered, Rudin's theorem implies that $C(\alpha K)^* \cong \ell^1(\alpha K)$. 
 
 Let $\iota: C_0(K) \to C(\alpha K)$ be the natural inclusion, so that $\iota^*: C(\alpha K)^* \to C_0(K)^*$ is an exact quotient map, meaning that it maps the closed unit ball onto the closed unit ball. Note that $\ker{\iota^*} = \spn \delta_\infty$, where $\delta_\infty$ is the evaluation functional at the point $\infty$. Thus, the First Isomorphism Theorem gives 
 \begin{equation*}
     C_0(K)^* \cong C(\alpha K)^* / \ker{\iota^*}  \cong \ell^1(\alpha K)/\spn {\delta_\infty}  \cong \ell^1(K).  \qedhere
 \end{equation*}
\end{proof}

\begin{definition}
    Let $\mu \in C_0(K)^*$, we call the subset 
    \begin{equation*}
        \supp({\mu}) = \{k \in K: \mu(\{k\}) \not = 0\}
    \end{equation*} 
    the \emph{support of $\mu$}.
\end{definition}

Observe that given $\mu \in C_0(K)^*$, we have $\supp({\mu}) \in \inset{K}$ and
\begin{equation*}
    \mu = \sum_{k \in \supp{(\mu)}} \mu(\{k\}) \delta_k.
\end{equation*}  

\begin{definition}
    For $C \subset K$, we define the \emph{restriction operator} $R_C:  C_0(K)^* \to C_0(K)^*$, acting on $\mu \in C_0(K)^*$ by
    \begin{equation*}
       R_C \mu = \sum_{k \in C} \mu(\{ k \})\delta_k.
    \end{equation*}
\end{definition}

It is easy to check that $R_C$ is a projection operator with $\norm{R_C} = 1$ (except in the trivial case $C = \emptyset$). If one thinks of the space $C_0(K)^* \cong \ell^1(K)$ as spanned by the transfinite Schauder basis $(\delta_k)_{k \in K}$, then the restriction operator $R_C$ is the standard basis projection. Therefore, we have the following properties.

\begin{lemma}\label{lmm: rest-disj-reg}
    Let $\{C_1, C_2, \dots, C_N\}$ be pairwise disjoint subsets of $K$ and set $C = \bigsqcup_{n=1}^N C_n$. Then:
    \begin{enumerate}[label={\normalfont{(\alph*)}}]
        \item  $R_{C} = \sum_{n=1}^N R_{C_n}$, in particular, $\norm{\sum_{n=1}^N R_{C_n}} = \norm{R_C} = 1$.
        \item For any $\mu \in C_0(K)^*$, we have $\norm{\sum_{n=1}^N R_{C_n} \mu} = \sum_{n=1}^N \norm{R_{C_n} \mu}$.
    \end{enumerate}
\end{lemma}

\bigskip
\section{Preliminary results}\label{sec-technical}

In this section, we introduce two ingredients that are key in the proof of Theorem \ref{th: main}. From now on, $X$ will always denote a Banach space and $T: C_0(K) \to X$ an operator with $\norm{T}_e = 1$. Furthermore, for the rest of this section, we fix an arbitrary $\delta \in (0, 1/2)$ and $\theta > 0$ small enough so that it satisfies $2(\delta + 2\theta) < 1$.

Before we move on to the first result, some notation is in order.

\begin{definition}
    Let $\eta \in C_0(K)$, we define the \emph{zero set of $\eta$} by 
    \begin{equation*}
        \mathbf{Z}(\eta) = \{k \in K: \eta(k) = 0 \} = \eta^{-1}(\{0\}),
    \end{equation*}
    and similarly 
    \begin{equation*}
        \nonvan{\eta} = K \backslash \mathbf{Z}(\eta) = \{k \in K: \eta(k) \not = 0 \}.
    \end{equation*}
    Further, we say that $\eta$ \emph{vanishes at a set $F \subset K$} if $F \subset \mathbf{Z}(\eta)$.
\end{definition}

Loosely speaking, our next result claims that we can find a sequence of functionals $(x^*_j)_{j=1}^\infty \subset B_{X^*}$ which can be decomposed in a way that makes them almost orthonormal and that witness the action of $T$.

These functionals will be the key ingredient in the proof of Theorem \ref{th: main}.

\begin{proposition}\label{prop-almost-orthonormal-system}
    There exist $(x^*_j)_{j=1}^\infty \subset B_{X^*}$ and finite sets $ F_j \subset \Tsupp{x^*_j}$ satisfying the following properties:
    \begin{enumerate}[label={\normalfont{(\alph*)}}]
        \item \label{prop-a} $F_i \cap F_j = \emptyset$ when $i \not = j$, or in other words, the family $(F_j)_{j=1}^\infty$ is pairwise disjoint.
        \item \label{prop-b} $\norm{R_{F_j} T^* x_j^*} > \delta + 2\theta$ for all $j \in \mathbb{N}$.
    \end{enumerate}
\end{proposition}

\begin{proof}
   Since $\norm{T^*} =\norm{T} \geq \norm{T}_e = 1$, there exists $x^*_1 \in B_{X^*}$ such that $\norm{T^* x^*_1} > \delta + 2\theta$. Choosing enough terms of the $\ell^1$ representation of $T^*x_1^*$ gives us the desired finite set $F_1 \subset \Tsupp{x^*_1}$.

   Assume we have built $(x^*_j)_{j=1}^J \subset B_{X^*}$ and finite sets $F_j \subset \Tsupp{x^*_j}$, such that $(F_j)_{j=1}^J$ is pairwise disjoint and \ref{prop-b} holds for $j=1, \dots, J$, we aim to build the next one.
    
   Let $A_{J} = \bigcup_{j=1}^J F_j \in \finset{K}$, and let $\mathcal{U}_J$ be a pairwise disjoint collection of neighbourhoods around $A_{J}$, which exists by Lemma \ref{lmm-neigh-collection}. Note that $\norm{T(I - \proj(\mathcal{U}_J, A_{J}))} \geq \norm{T}_e = 1$ by Lemma \ref{lmm: proj} \ref{lmm: proj-a}. Thus there exists $\zeta_{J+1} \in B_{C_0(K)}$ such that $\norm{T(I - \proj(\mathcal{U}_J, A_{J}))\zeta_{J+1}} > 2(\delta + 2\theta)$. 
    
   Define $\eta_{J+1} = \frac{1}{2}(I-\proj(\mathcal{U}_J, A_{J}))\zeta_{J+1}$ so that $\norm{T \eta_{J+1}} > \delta + 2\theta$.  By the Hahn-Banach theorem, there exists $x^*_{J+1} \in B_{X^*}$ such that $x^*_{J+1} T \eta_{J+1} > \delta + 2\theta$.

    Since $\eta_{J+1} \in B_{C_0(K)}$ by Lemma \ref{lmm: proj} \ref{lmm: proj-d}, it follows that 
   \begin{align*}
       \norm{R_{\nonvan{\eta_{J+1}}} T^*x^*_{J+1} } &\geq \lvert \langle R_{\nonvan{\eta_{J+1}}} T^*x^*_{J+1}, \eta_{J+1} \rangle \rvert
        =  \lvert \langle T^*x^*_{J+1}, \eta_{J+1} \rangle \rvert \\ 
        &> \delta + 2\theta.
   \end{align*}
   In particular, there exists a finite set of points $F_{J+1} \subset \Tsupp{x^*_{J+1}} \cap \nonvan{\eta_{J+1}}$ such that $\norm{R_{F_{J+1}} T^*x^*_{J+1}} > \delta + 2 \theta$. Finally, by Lemma \ref{lmm: proj} \ref{lmm: proj-b} we have $\nonvan{\eta_{J+1}} \cap A_{J} = \emptyset$, so $F_{J+1}$ is disjoint with $F_{1}, \dots, F_{J}$. This finishes the recursive construction.
\end{proof}

\begin{rems}\label{rem: factor-2}
    For later reference, we observe that the condition $\delta < 1/2$ stems from the fact that, to define $\eta_{J+1}$, we need to consider $(I - \proj(\mathcal{U}_J, A_J))\zeta_{J+1}$ whose norm may be greater than one, but it is at most two by Lemma \ref{lmm: proj} \ref{lmm: proj-d}.
\end{rems}

Now, for each $n \in \mathbb{N}$, we aim to separate the points $F_n$ from the points in $\bigcup_{j \not = n} F_j $ by a collection of neighbourhoods. Once we do that, we are set, since we can take the sum of the indicator functions of those neighbourhoods. However, some technical care is required. \\

Roughly speaking, the idea is that, if we cannot (almost) separate the points in $F_n$ from those in $\bigcup_{j \not = n} F_j $, it is because a significant amount of weight is accumulating around the points in $F_n$. However, since the total weight is finite, this can only happen for a finite number of sets $F_n$'s, thus, removing those if necessary, we can (almost) separate the remaining sets.

\begin{lemma}\label{lmm-ess-weight}
    Let $(x^*_j)_{j=1}^\infty \subset B_{X^*}$ and let $(F_j)_{j=1}^\infty \subset \finset{K}$ be pairwise disjoint.  Then, for any $\varepsilon > 0$, there exist $n \in \mathbb{N}$ and a pairwise disjoint neighbourhood collection $\mathcal{U}_n$ around $F_n$ with the property that for all $J \in \mathbb{N}$,  there exists $j > J$ with $\norm{R_{S(\mathcal{U}_n)} T^* x^*_j} < \varepsilon$.
\end{lemma}

\begin{proof}
    We proceed by contradiction and assume that the claim is false for some $\varepsilon > 0$ and choose $N > \norm{T}/\varepsilon$.

    Let $F = \bigsqcup_{n=1}^N F_n$. By Lemma \ref{lmm-neigh-collection} there exists a pairwise disjoint neighbourhood collection $\mathcal{U}$ around $F$, which naturally yields a pairwise disjoint neighbourhood collection $\mathcal{U}_n$ around $F_n$ for each $n=1, \dots, N$, with the property $\mathcal{U} = \bigcup_{n=1}^N \mathcal{U}_n$.

    Since the claim is false, for each $n=1, \dots, N$, there exists $J_n$ such that $\norm{R_{S(\mathcal{U}_n)} T^* x^*_j} \geq \varepsilon$ for all $j \geq J_n$. 

    Let $J = \max \{J_1, \dots, J_N\}$, so that for $x^*_{J}$ we have
    \begin{equation*}
        \norm{T} \geq \norm{T^* x^*_{J}} \geq \norm{\sum_{n=1}^N R_{S(\mathcal{U}_n)} T^* x^*_{J}} = \sum_{n=1}^N \norm{ R_{S(\mathcal{U}_n)} T^* x^*_{J}} \geq N \varepsilon,
    \end{equation*}
    by using Lemma \ref{lmm: rest-disj-reg}. This contradicts our choice of $N$, thus proving the proposition. 
\end{proof}

\begin{corollary}\label{corollary-ess-weight}
     Let $(x^*_j)_{j=1}^\infty \subset B_{X^*}$ and let $(F_j)_{j=1}^\infty \subset \finset{K}$ be pairwise disjoint. Then for any $\varepsilon > 0$, there exist $n \in \mathbb{N}$, a pairwise disjoint neighbourhood collection $\mathcal{U}_n$ around $F_n$ and a subsequence $(x^*_{k_j})_{j=1}^\infty$ of $(x^*_{j})_{j=1}^\infty$ such that $\lVert R_{S(\mathcal{U}_n)} T^* x^*_{k_j} \rVert < \varepsilon$ for all $j \in \mathbb{N}$. 
\end{corollary}
\bigskip
\section{Proof of Theorem \ref{th: main}}\label{sec-main-proof}

Informally, we aim to find a sequence $(A_j)_{j=1}^\infty$ of finite subsets of $K$  which can be separated by pairwise disjoint neighbourhoods collections $(\mathcal{V}_j)_{j=1}^\infty$, such that they witness the action of $T^* y^*_j$ for some functionals $y^*_j$, in a way that will be made precise in Lemma \ref{lmm: auxilary-main}.

Using this, we can build the functions $(\eta_j)_{j=1}^\infty$ to be supported in the neighbourhoods belonging to $\mathcal{V}_j$ in a natural way. The fact that $A_j$ witnesses the action of $T^* y^*_j$ will guarantee that $\norm{T \eta_j}  > \delta$.

We now proceed to the construction of the family of points and neighbourhoods, by using a diagonal argument on the family given by Proposition \ref{prop-almost-orthonormal-system}.

\begin{lemma}\label{lmm: auxilary-main}
    There exist $(y^*_j)_{j=1}^\infty \subset B_{X^*}$, finite sets $A_j \subset \Tsupp{y^*_j}$, and pairwise disjoint neighbourhood collections $\mathcal{V}_j$ around $A_j$ satisfying:
    
    \begin{enumerate}[label={\normalfont{(\alph*)}}]
        \item \label{lmm: disj-A} The family $\bigcup_{j=1}^\infty \mathcal{V}_j$ is pairwise disjoint; that is $U \cap V = \emptyset$ for any distinct $U, V \in \bigcup_{j=1}^\infty \mathcal{V}_j$.
        \item \label{lmm: pesos-delta} $\norm{R_{A_j} T^* y_j^*} > \delta + \theta$ for each $j \in \mathbb{N}$.
    \end{enumerate}
\end{lemma}

\begin{proof}
    By recursion, for each $J  \in \mathbb{N}_0$, we shall construct
    a sequence of functionals $(x^*_{i,J})_{i=1}^\infty$, with associated finite sets $F_{i,J} \subset \Tsupp{x^*_{i,J}}$, a functional $y^*_J \in B_{X^*}$, a finite set $A_J \subset \Tsupp{y^*_J}$ and a neighbourhood collection $\mathcal{V}_J$ around $A_J$, satisfying:
    \begin{enumerate}[label={\normalfont{(\roman*)}}]
        \item\label{it: cond1} $\norm{R_{F_{i,J}} T^* x^*_{i,J}} > \delta + 2\theta - \sum_{n=1}^J \theta/2^n$ for every $i \in \mathbb{N}$ and the family $(F_{i,J})_{i=1}^\infty$ is pairwise disjoint.
        \item\label{it: cond2}  $\norm{R_{A_{J}} T^* y^*_{J}} > \delta + \theta$ provided that $J \geq 1$ (thus, the condition becomes trivial in the base case $J=0$).
        \item\label{it: cond3} The family $\bigcup_{j=0}^J \mathcal{V}_j$ is pairwise disjoint and $F_{i,J} \cap \left(\bigcup_{j=0}^J \mathcal{V}_j\right) = \emptyset$ for all $i \in \mathbb{N}$.
    \end{enumerate}
    By Proposition \ref{prop-almost-orthonormal-system}, there exists a sequence of functionals $(x^*_{i,0})_{i=1}^\infty$, with associated finite sets $F_{i,0} \subset \Tsupp{x^*_{i,0}}$ satisfying \ref{it: cond1}. Taking $y^*_0 = 0$, $A_0 = \emptyset$, and $\mathcal{V}_0 = \emptyset$ makes \ref{it: cond2} and \ref{it: cond3} trivially hold.

    Suppose we have carried out the construction for $j = 0, 1, \dots, J$, and we will show that we can continue it to $J+1$. By Corollary \ref{corollary-ess-weight} applied to $(x^*_{i,J})_{i=1}^\infty$ and $(F_{i,J})_{i=1}^\infty$, there exist $n = n_{J+1} \in \mathbb{N}$, a pairwise disjoint neighbourhood collection $\mathcal{U}_{n}$ around $F_{n,J}$ as well as a subsequence $(x^*_{k_i,J})_{i=1}^\infty$ of  $(x^*_{i,J})_{i=1}^\infty$ such that $\lVert R_{S(\mathcal{U}_{n})} T^* x^*_{k_i,J} \rVert < \theta/2^{J+1}$ for all $i \in \mathbb{N}$.

    By shrinking the neighbourhoods in $\mathcal{U}_{n}$, we can assume that $\mathcal{U}_{n} \cup \left(\bigcup_{j=0}^J \mathcal{V}_j\right)$ is pairwise disjoint; note that this does not affect the condition that $\lVert R_{S(\mathcal{U}_{n})} T^* x^*_{k_i,J} \rVert < \theta/2^{J+1}$ for all $i \in \mathbb{N}$. 

    Define $y^*_{J+1} = x^*_{n,J}$, $A_{J+1} = F_{n,J}$, $\mathcal{V}_{J+1} = \mathcal{U}_{n}$, as well as sequences $x^*_{i, J+1} = x^*_{k_i, J}$ and $F_{i, J+1} = F_{{k_i}, J} \backslash S(\mathcal{U}_{n})$ for $i \in \mathbb{N}$.
    
    Then conditions \ref{it: cond1}, \ref{it: cond2} and \ref{it: cond3} are satisfied for $J+1$, which proves that we can carry on the recursive construction and thus the result follows.
\end{proof}

As already mentioned, by using the sequence $(A_j)_{j=1}^\infty$ and the neighbourhood collections $(\mathcal{V}_j)_{j=1}^\infty$ from the previous lemma, it is easy to build a sequence of disjoint functions $(\eta_j)_{j=1}^\infty$, each of them supported in $\mathcal{V}_j$ in a natural way. For technical reasons, we may need to further shrink the neighbourhoods in $\mathcal{V}_j$ to be able to control the image of the function $\eta_j$ under $T$. \\

The reader may wish to compare the following proof with that of Theorem \ref{th: pelc}, and observe that the additional assumption that $K$ is scattered makes it possible not only to explicitly construct the family $(\eta_j)_{j=1}^\infty$, but also to obtain a lower bound for the norms of their images in terms of the essential norm of the operator.

\begin{proof}[Proof of Theorem \ref{th: main}]

Let $(y^*_j)_{j=1}^\infty$, $(A_j)_{j=1}^\infty$ and $(\mathcal{V}_j)_{j=1}^\infty$ be as in Lemma \ref{lmm: auxilary-main}. For each $j \in \mathbb{N}$, partition the set $\Tsupp{y^*_j} = A_j \sqcup H_j \sqcup G_j$, in such a way that $H_j \in \finset{K}$ and $\norm{R_{G_j} T^* y^*_j} < \theta$. Further, by shrinking the open neighbourhoods in the family $\mathcal{V}_j$ if necessary, we may assume that $S(\mathcal{V}_j) \cap H_j = \emptyset$. Note that this does not change the pairwise disjointness of the family $\bigcup_{j=1}^\infty \mathcal{V}_j$.

For each $k \in A_j$, let $\lambda_k$ be an unimodular scalar such that the equation $\lambda_k T^*y^*_{j} (\{k\}) = \lvert T^*y^*_{j} (\{k\}) \rvert$  holds and define $\eta_j = \sum_{k \in A_{j}} \lambda_k \id{V_j(k)}$, where $V_j: A_j \to \mathcal{V}_j$ is the bijection given in Definition \ref{def: nhg-collection}.

The pairwise disjointness of $\mathcal{V}_j$ guarantees that $\norm{\eta_j} \leq 1$, while the pairwise disjointness of $\bigcup_{j=1}^\infty \mathcal{V}_j$ guarantees that the functions $(\eta_j)_{j=1}^\infty$ are disjoint. 

Further, we have
\begin{align*}
    \langle y^*_{j}, T\eta_j \rangle  &= \sum_{k \in A_{j}} \lambda_k \langle T^* y^*_{j}, \id{V_j(k)} \rangle  = \sum_{k \in A_{j}} \sum_{h \in \Tsupp{y^*_{j}}}  \lambda_k T^* y^*_{j} (\{h\}) \id{V_j(k)}(h).
\end{align*}
We can split the sum into three terms, by using the partition $\Tsupp{y^*_{j}} = A_{j} \sqcup H_{j} \sqcup G_{j}$. For the first term, note that
\begin{align*}
    \sum_{k \in A_{j}} \sum_{h \in A_{j}}  \lambda_k T^* y^*_{j} (\{h\}) \id{V_j(k)}(h) = \sum_{k \in A_{j}} \lvert T^* y^*_{j} (\{k\}) \rvert = \norm{R_{A_{j}} T^* y^*_{j}} > \delta + \theta.
\end{align*}
The second term is zero since $\id{V_j(k)}(h) = 0$ when $h \in H_{j}$ and $k \in A_j$ because $S(\mathcal{V}_j) \cap H_j = \emptyset$. Finally, for the third term, we have
\begin{align*}
    \left|\sum_{k \in A_{j}} \sum_{h \in G_{j}}  \lambda_k T^* y^*_{j} (\{h\}) \id{V_j(k)}(h) \right| \leq \norm{R_{G_{j}} T^* y^*_{j}} < \theta,
\end{align*}
where we used that for each $h \in G_j$ there is at most one $k \in A_j$ such that $h \in V_j(k)$, by the pairwise disjointness of $\mathcal{V}_j$. Overall, we get that
\begin{equation*}
    \norm{T\eta_j} \geq |\langle y^*_{j}, T\eta_j \rangle| > (\delta + \theta) - \theta = \delta,
\end{equation*}
which completes the proof. 
\end{proof} 
\bigskip
\section{Factorizations and the essential norm of the adjoint}\label{sec-sharpness and improvements}

It turns out that the natural way of measuring how well the identity on $c_0$ factors through an operator $T$ has to do with the essential norm of the adjoint $T^*$, rather than with the essential norm of $T$ itself. As discussed in Remark \ref{rem: factor-2}, the condition $\delta < 1/2$ stems from Proposition \ref{prop-almost-orthonormal-system}.

Therefore, it suffices to establish an analogue of this proposition where $\delta$ can be taken in the interval $(0, \norm{T^*}_e)$. Once we have done so, the rest of the proof for the factorization carries over.

\begin{proposition}\label{analogous-prop-orth}
    Let $\delta \in (0, \norm{T^*}_e)$ and fix $\theta$ such that $\delta + 2\theta < \norm{T^*}_e$. Then there exist $(x^*_j)_{j=1}^\infty \subset B_{X^*}$ and finite sets $ F_j \subset \Tsupp{x^*_j}$ satisfying conditions \ref{prop-a} and \ref{prop-b} of Proposition \ref{prop-almost-orthonormal-system}.
\end{proposition}
\begin{proof}
   We carry out the construction by recursion as in the proof of Proposition \ref{prop-almost-orthonormal-system}. As before, we can find $x_1^* \in B_{X^*}$ and a finite $F_1 \subset \Tsupp{x^*_1}$ such that $\norm{T^* x_1^*} > \delta + 2\theta$. 
   Since $\norm{(I - R_{F_1})T^*} \geq \norm{T^*}_e > \delta + 2\theta$, we can find $x^*_2 \in B_{X^*}$ such that $\norm{(I - R_{F_1})T^*x_2^*} > \delta + 2\theta$, so that there exists a finite $F_2 \subset \Tsupp{x^*_2}$ such that $F_2 \cap F_1 = \emptyset$ and $\norm{R_{F_2}T^* x^*_2} > \delta + 2\theta$. Continuing in this way finishes the construction.
\end{proof}

\begin{rems}
    The essence of the argument is the following: while in $C_0(K)$ the possible presence of accumulation points only allows for the bound $\norm{I - \proj(\mathcal{U}, F)} \leq 2$, in the dual case $C_0(K)^* \cong \ell^1(K)$ we have a natural way to project. Namely, in the dual formulation, the role of $\proj(\mathcal{U}, F)$ is played by $R_{F}$ and thus we can guarantee $\norm{I - R_F} = 1$ by Lemma \ref{lmm: rest-disj-reg}.
\end{rems}

We can now invoke Proposition \ref{analogous-prop-orth} in place of Proposition \ref{prop-almost-orthonormal-system} when proving Lemma \ref{lmm: auxilary-main}. This allows us to establish Lemma \ref{lmm: auxilary-main} under the condition $\delta + \theta < \norm{T^*}_e$, instead of the previous requirement $\delta + \theta < \norm{T}_e/2$. Consequently, the proof of Theorem \ref{th: main} carries through for any $\delta < \norm{T^*}_e$, which, when combined with Theorem \ref{thr-niels}, yields the following.

\begin{corollary}
    Let $K$ be a scattered, locally compact Hausdorff space and let $X$ be a Banach space that does not contain a copy of $\ell^1$ or for which the unit ball of $X^*$ is weak$^*$ sequentially compact. Then for any non-compact operator $T: C_0(K) \to X$ and any $\varepsilon > 0$, there exist operators $U \in \mathscr{B}(X, c_0)$ and $V \in \mathscr{B}(c_0, C_0(K))$ such that
    \begin{equation*}
        UTV = I_{c_0} \hspace{5pt} \text{ and } \hspace{5pt} \norm{U} \norm{V} < \frac{1}{\norm{T^*}_e} + \varepsilon.
    \end{equation*}
\end{corollary}

Now we can prove Theorem \ref{th: ess-norm-adjoint}.

\begin{proof}[Proof of Theorem \ref{th: ess-norm-adjoint}]
    Our previous corollary shows that for any $\varepsilon > 0$, there exist $U \in \mathscr{B}(X, c_0)$ and $ V \in \mathscr{B}(c_0, C_0(K))$ such that $UTV = I_{c_0}$ and $\norm{U}\norm{V} < 1/{\norm{T^*}_e} + \varepsilon$, in other words, $\norm{U}^{-1}\norm{V}^{-1} > \norm{T^*}_e/{(\varepsilon \norm{T^*}_e + 1)}$. Since this is true for any $\varepsilon > 0$, we get
    \begin{equation*}
        \norm{T^*}_e \leq \sup \{\norm{U}^{-1}\norm{V}^{-1}: UTV = I_{c_0} \},
    \end{equation*}
    and we are left to show the reverse inequality.

    Suppose now that $U \in \mathscr{B}(X, c_0)$ and $V \in \mathscr{B}(c_0, C_0(K))$ satisfy $UTV = I_{c_0}$. Then for any $S \in \mathscr{K}(X^*, C_0(K)^*)$, we get
    \begin{align*}
         \norm{V^*}\norm{T^* - S}\norm{U^*} &\geq \norm{V^*(T^* - S)U^*} = \norm{I_{c_0}^* - V^*S U^*}
         \geq \norm{I_{c_0}^*}_e = 1,
    \end{align*}
    so it follows that $\norm{T^* - S} \geq \norm{U}^{-1} \norm{V}^{-1}$. Taking the infimum over all compact operators gives $\norm{T^*}_e \geq \norm{U}^{-1} \norm{V}^{-1}$, which proves the upper bound in the supremum.
\end{proof}

\bigskip
\section{Minimality of the essential norm }\label{sec-minimality}

This section is dedicated to the proof of Theorems \ref{th: main-minimality} and \ref{th: main-bonsall-min}, which will follow as a consequence of our factorization results in the form of Theorems \ref{th: fact-iden} and \ref{th: ess-norm-adjoint}, respectively. We start with the following known result.

\begin{lemma}\label{lmm: complemented-c0}
    Let $K$ be a scattered, locally compact Hausdorff space. Then, there exist $A \in \mathscr{B}(c_0, C_0(K))$ and $B \in \mathscr{B}(C_0(K), c_0)$ such that
    \begin{equation*}
        BA = I_{c_0} \hspace{10pt} \text{ and } \norm{A} = \norm{B} = 1.
    \end{equation*}
\end{lemma}
\begin{proof}
    Recall that the unit ball of $C_0(K)^*$ is weak$^*$ sequentially compact (see \cite[Theorem 6.2]{arnott2023uniqueness}), while $C_0(K)$ always contains an isometric copy of $c_0$. The result now follows from \cite[Theorem 6]{dowling1999remarks}.
\end{proof}

We will also need the following technical lemma.

\begin{lemma}\label{lmm: multiple-factori}
    There exists a sequence of non-compact operators $(P_n)_{n=1}^\infty \subset \mathscr{B}(C_0(K))$ such that, for every sequence of operators $(T_n)_{n=1}^\infty \subset \mathscr{B}(C_0(K))$ with $\norm{T_n}_e > 2^{2n+1}$ for every $n \in \mathbb{N}$, there exist $U, V \in \mathscr{B}(C_0(K))$ such that 
    \begin{equation}\label{eq: sim-factorization}
        P_n U T_n V P_n = P_n \hspace{10 pt} \text{ for all } n \in \mathbb{N}.
    \end{equation}
\end{lemma}
    \begin{proof}
    
        We note that $c_0 \cong \left(\bigoplus_{n=1}^\infty c_0\right)_{c_0} $ and let
    \begin{equation*}
        \pi_n: \left(\bigoplus_{n=1}^\infty c_0\right)_{c_0}  \to c_0 \hspace{5pt} \text{ and } \hspace{5pt} \iota_n: c_0 \to \left(\bigoplus_{n=1}^\infty c_0\right)_{c_0} 
    \end{equation*}  
    be the projection onto the $n$th term of the sum and the inclusion into the $n$th term, so that $\pi_n \circ \iota_n = I_{c_0}$ and $\norm{\pi_n} = \norm{\iota_n} = 1$.

    From Lemma \ref{lmm: complemented-c0} it follows in particular that $c_0$ is complemented in $C_0(K)$. Thus, there exists a projection \break $\pi \in \mathscr{B}\left(C_0(K), \left(\bigoplus_{n=1}^\infty c_0\right)_{c_0} \right)$ of $C_0(K)$ onto $\left(\bigoplus_{n=1}^\infty c_0\right)_{c_0}$ and an inclusion $\iota \in \mathscr{B}\left(\left(\bigoplus_{n=1}^\infty c_0\right)_{c_0}, C_0(K)\right)$ of $\left(\bigoplus_{n=1}^\infty c_0\right)_{c_0}$ into $C_0(K)$, so that the composition $\pi \circ \iota$ gives the identity map on $\left(\bigoplus_{n=1}^\infty c_0\right)_{c_0}$.
    
    Let $P_n = \iota \circ \iota_n \circ \pi_n \circ \pi \in \mathscr{B}(C_0(K))$ for each $n \in \mathbb{N}$. Since $C_0(K)$ has weak$^*$ sequentially compact dual unit ball, we can apply Theorem \ref{th: fact-iden} to each $T_n$. Thus, for each $n \in \mathbb{N}$ there exist $U_n \in \mathscr{B}(C_0(K), c_0)$ and $V_n \in \mathscr{B}(c_0, C_0(K))$ such that
    \begin{equation*}
        U_n T_n V_n = I_{c_0} \hspace{5pt} \text{ and } \hspace{5pt}  \norm{U_n}\norm{V_n} < 1/2^{2n},
    \end{equation*}
    and by re-scaling we may assume both $\norm{U_n} < 1/2^n$ and $\norm{V_n} < 1/2^n$.

    Note that $V_n \pi_n \pi \in \mathscr{B}(C_0(K))$ and $\norm{V_n \pi_n \pi} \leq \norm{\pi} \norm{V_n} < \norm{\pi}/2^n$, so that we can define the operator $V = \sum_{n=1}^\infty V_n \pi_n \pi$.

    Similarly, $\iota \iota_n U_n \in \mathscr{B}(C_0(K))$ and $\norm{\iota \iota_n U_n } \leq \norm{\iota}\norm{U_n} < \norm{\iota}/2^n$, so again we can define the operator $U = \sum_{n=1}^\infty \iota \iota_n U_n $.

    From the definition, it is straightforward to check that we have $\pi_m \circ \iota_n = \delta_{n,m} I_{c_0}$ and thus elementary computations show $P_n U = \iota \iota_n U_n$ and $V P_n = V_n \pi_n \pi$, which gives
    \begin{equation*}
        P_n U T_n V_n P_n = \iota \iota_n U_n T_n V_n \pi_n \pi = \iota \iota_n \pi_n \pi = P_n.
    \end{equation*}
    
    \begin{figure}[H]
    \centering
    \includegraphics[scale = 0.95]{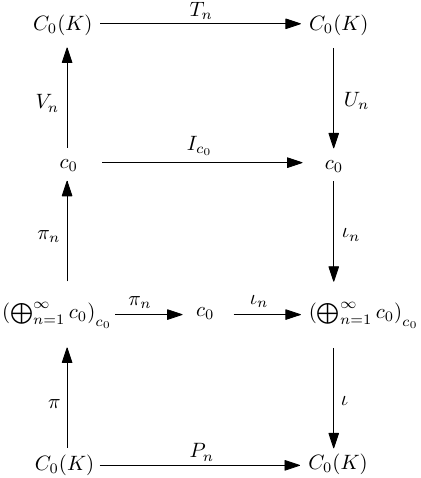}
    \end{figure}

\end{proof}

In the next proofs, given a norm $\vertiii{\cdot}$ on $\mathscr{B}(C_0(K))/\mathscr{K}(C_0(K))$, we refer indiscriminately to $\vertiii{T}$ of any $T \in \mathscr{B}(C_0(K))$. This is to be understood as the seminorm this norm induces on $\mathscr{B}(C_0(K))$ or, alternatively, as the norm of the equivalence class represented by $T$. Using this convention and equipped with the previous lemma, we can prove the minimality of the Calkin algebra norm. 

\begin{proof}[Proof of Theorem \ref{th: main-minimality}]
    Choose $(P_n)_{n=1}^\infty$ as in the proof of Lemma \ref{lmm: multiple-factori}. We proceed by contradiction and assume there exists a norm $\vertiii{\cdot}$ on the Calkin algebra such that for any $\varepsilon > 0$, we can find $R \in \mathscr{B}(C_0(K))$ such that $\norm{R}_e = 1$ and $\vertiii{R} < \varepsilon$.

    Thus, for each $n \in \mathbb{N}$, there exists an operator $R_n \in \mathscr{B}(C_0(K))$ with $\norm{R_n}_e = 1$ and
    \begin{equation*}
         \vertiii{R_n} <  \frac{1}{n (2^{2n+1} + 1) \vertiii{P_n}},
    \end{equation*} 
    where we emphasize that $\vertiii{P_n} \not = 0$ since $P_n \not \in \mathscr{K}(C_0(K))$. Define $T_n = (2^{2n+1} + 1)R_n$, so that $\norm{T_n}_e > 2^{2n+1}$.
    
    By Lemma \ref{lmm: multiple-factori} applied to the family $(T_n)_{n=1}^\infty$, there exists $U,V \in \mathscr{B}(C_0(K))$ such that
    \begin{equation*}
        P_n U T_n V P_n = P_n \hspace{10pt} \text{ for all } n \in \mathbb{N},
    \end{equation*}
    and thus
    \begin{equation*}
        \vertiii{P_n} = \vertiii{P_n U T_n V P_n} \leq \vertiii{P_n}^2 \vertiii{U} \vertiii{V} \vertiii{T_n}.
    \end{equation*}
    It follows that
    \begin{equation*}
        0 < \frac{1}{\vertiii{U}\vertiii{V}} \leq \vertiii{P_n} \vertiii{T_n} = (2^{2n + 1} + 1)\vertiii{P_n} \vertiii{R_n} < \frac{1}{n}, 
    \end{equation*}
    and making $n \to \infty$ yields a contradiction, finishing the proof.
\end{proof}

\begin{rems} 
The previous argument works in a more general setting. Namely, let $X$, $Y$ be Banach spaces. For $T \in \mathscr{B}(X)$ define
\begin{equation*}
    \nu_Y{(T)} = \sup\{ \norm{U}^{-1}\norm{V}^{-1}: UTV = I_Y\}
\end{equation*}
if $I_Y$ factors through $T$ and $\nu_Y{(T)} = 0$ otherwise.  Further, let $\mathcal{J}_Y = \{T \in \mathscr{B}(X): \nu_Y(T) = 0\}$, in other words, the set of operators that do not factor $I_Y$. Observe that the set $\mathcal{J}_Y$ is always closed.

Suppose that $\mathcal{J}_Y$ is closed under addition (and thus, a closed ideal of $\mathscr{B}(X)$) and that $Y$ admits a Schauder decomposition $(Y_n)_{n=0}^\infty$ with $\sup \{d_{BM}(Y_n, Y) : n \geq 1\} < \infty$, where $d_{BM}$ denotes the Banach-Mazur distance. Then $\nu_Y$ minorizes any algebra norm on $\mathscr{B}(X)/\mathcal{J}_Y$. In particular, if $\nu_Y$ defines an algebra norm, it is minimal.

Further, note that the conditions on $Y$ are satisfied, for example, by the classical sequence spaces $\ell^p$, $1 \leq p \leq \infty$.
\end{rems}

\begin{rem}
    These factorisation-type norms arise naturally in algebras of the form $\mathscr{B}(X)/\mathscr{I}$ for some ideal $\mathscr{I}$ and play a central role in the study of their norm properties. For example, in a forthcoming paper, Johnson and Phillips \cite{johnsonandphillips} introduce the notion of \emph{uniform incompressibility} (see also \cite[Definition 2.8]{arnott2023uniqueness} for a definition), and employ factorization-type norms to show this property holds in multiple Banach algebra of the form $\mathscr{B}(X)/\mathscr{I}$. Our results show that the essential norm is uniformly incompressible for the Calkin algebra of $C_0(K)$ spaces, $K$ a scattered, locally compact Hausdorff space. 
\end{rem}

Lastly, we provide a proof for Theorem \ref{th: main-bonsall-min}.

\begin{proof}[Proof of Theorem \ref{th: main-bonsall-min}]
    Let $\vertiii{\cdot}$ be any norm on $\mathscr{B}(C_0(K)) / \mathscr{K}(C_0(K))$ satisfying $\vertiii{S} \leq \norm{S^*}_e$ for all $S \in \mathscr{B}(C_0(K))$, we will show that then $\vertiii{S} = \norm{S^*}_e$ for all $S \in \mathscr{B}(C_0(K))$. For this, it is enough to show that $\vertiii{S} \geq \norm{S^*}_e$  for all $S \in \mathscr{B}(C_0(K))$.

    So fix any non-compact operator $T \in \mathscr{B}(C_0(K))$ and assume without loss of generality that $\norm{T^*}_e = 1$, we need to show that $\vertiii{T} \geq 1$.
    
    Choose $A \in \mathscr{B}(c_0, C_0(K))$ and $B \in \mathscr{B}(C_0(K), c_0)$ according to Lemma \ref{lmm: complemented-c0} and fix any $\varepsilon > 0$. By Theorem \ref{th: ess-norm-adjoint} we can find $V \in \mathscr{B}(c_0, C_0(K))$ and $U \in \mathscr{B}(C_0(K), c_0)$ such that $\norm{U} \norm{V} < 1 + \varepsilon $ and $UTV = I_{c_0}$, which gives 
    \begin{equation*}
        AUTVB = AB.
    \end{equation*}
   Note that since $I_{c_0} = BA$ then $(AB)^2 = AB \not = 0$, and thus
    \begin{equation*}
        1 \leq \vertiii{AB} = \vertiii{AUTVB} \leq \vertiii{AU}\vertiii{T}\vertiii{VB} \leq \norm{(AU)^*}_e \norm{(VB)^*}_e \vertiii{T} , 
    \end{equation*}
    where the last inequality follows since $\vertiii{S} \leq \norm{S^*}_e$ for all $S \in \mathscr{B}(C_0(K))$. Since
    \begin{align*}
        \norm{(AU)^*}_e \norm{(VB)^*}_e &\leq \norm{(AU)^*} \norm{(VB)^*} \leq \norm{A}\norm{U}\norm{V}\norm{B}
        \leq (1 + \varepsilon),
    \end{align*}
    we get
    \begin{equation*}
       1 \leq \norm{(AU)^*}_e \norm{(VB)^*}_e \vertiii{T} \leq (1 + \varepsilon) \vertiii{T}.
    \end{equation*}
    As this is true for any $\varepsilon > 0$, it follows that $1 \leq \vertiii{T}$, which proves the result.
\end{proof}

\bigskip
\section{Automatic continuity of homomorphism from $\mathscr{B}(C([0,\alpha]))$}\label{sec-uniqueness-norm}

We show how a simple argument allows us to extend Ogden's result on automatic continuity to all ordinals. We first establish some auxiliary results.

For a Banach algebra $\mathcal{B}$, we denote by $\mathbb{M}_n(\mathcal{B})$ the Banach algebra of $(n \times n)$-matrices whose entries are elements of $\mathcal{B}$.
We write $T = [T_{i,j}]_{i,j=1}^n  \in \mathbb{M}_n(\mathcal{B})$ for the matrix having $T_{i,j} \in \mathcal{B}$ for its $(i,j)$-entry. Denote by $\pi_{i,j}: \mathbb{M}_n(\mathcal{B}) \to \mathcal{B}$, $[T_{r,s}]_{r,s=1}^n \mapsto T_{i,j}$ (the map sending a matrix to its $(i,j)$-entry) and $\iota_{i,j}: \mathcal{B} \to \mathbb{M}_n(\mathcal{B})$ the map sending an element $T \in \mathcal{B}$ to the matrix with all entries zero except $T$ in the $(i,j)$-entry. We assume $\mathbb{M}_n(\mathcal{B})$ is equipped with any norm $\norm{\cdot}$ satisfying
\begin{equation*}
    \max \{ \norm{T_{i,j}}: 1 \leq i,j \leq n\} \leq \norm{T} \leq \sum_{i,j = 1}^n \norm{T_{i,j}},
\end{equation*}
so that both $\iota_{i,j}$ and $\pi_{i,j}$ are continuous.

\begin{proposition}
    Let $\mathcal{B}$ be a unital Banach algebra. Then every homomorphism from $\mathcal{B}$ to a Banach algebra is continuous if and only if every homomorphism from $\mathbb{M}_n(\mathcal{B})$ to a Banach algebra is continuous for every $n \geq 1$.
\end{proposition}
\begin{proof}
    The reverse direction is immediate, so we only need to show that if every homomorphism from $\mathcal{B}$ to a Banach algebra is continuous, the same holds for $\mathbb{M}_n(\mathcal{B})$, $n \geq 2$.
    
    Let $\mathcal{A}$ be a Banach algebra and $\theta: \mathbb{M}_n(\mathcal{B}) \to \mathcal{A}$ be a homomorphism. Since $\sum_{i,j=1}^n \iota_{i,j} \circ \pi_{i,j} = \identity_{\mathbb{M}_n(\mathcal{B})}$, it follows that
    \begin{equation*}
        \theta = \sum_{i,j = 1}^n (\theta \circ \iota_{i,j}) \circ \pi_{i,j},
    \end{equation*}
    and since $\pi_{i,j}$ is continuous, to prove continuity of $\theta$, it is enough to show continuity of $\theta_{i,j} = \theta \circ \iota_{i,j}: \mathcal{B} \to \mathcal{A}$ for each pair $1 \leq i,j \leq n$.
    
    Note that $\theta_{i,j}$ is linear. Furthermore, for $1 \leq i \leq n$, $\iota_{i,i}$ is an algebra homomorphism; therefore $\theta_{i,i}$ is an algebra homomorphism and thus continuous by hypothesis.

    We focus now on the off-diagonal terms. Let $i \not = j$ be fixed. For any $S, U \in \mathcal{B}$ we have
    \begin{align*}
        \theta_{i,j}(SU) &= \theta(\iota_{i,j}(SU)) = \theta(\iota_{i,i}(S)\iota_{i,j}(U)) \\
        &= \theta(\iota_{i,i}(S)) \theta(\iota_{i,j}(U)) = \theta_{i,i}(S) \theta_{i,j}(U).
    \end{align*}
    For any null-sequence $(S_k)_{k \in \mathbb{N}}$, it follows that
    \begin{equation*}
        \theta_{i,j}(S_k) = \theta_{i,i}(S_k) \theta_{i,j}(\identity_\mathcal{B}) \to 0
    \end{equation*}
    as $k \to \infty$, where we used that $\theta_{i,i}$ is continuous. This proves the continuity of $\theta_{i,j}$ and finishes the proof.
\end{proof}

We remark in passing that the assumption that the Banach algebra $\mathcal{B}$ is unital can be substantially weakened. In particular, it is enough
that null sequences in $\mathcal{B}$ factor.

Expressing the operators on $X^n = X \oplus X \oplus \dots \oplus X$ in matrix form, we have $\mathscr{B}(X^n) = \mathbb{M}_n(\mathscr{B}(X))$ and thus we get the following.

\begin{corollary}\label{cor: maximality norm}
    Let $X$ be a Banach space. Then the operator norm is maximal in $\mathscr{B}(X)$ if and only if the operator norm is maximal in $\mathscr{B}(X^n)$ for every $n \geq 1$.
\end{corollary}

We also recall the classification theorem for spaces of continuous functions on ordinals, originally proven by \cite{gul1975isomorphic, kislyakov1975classification}; the specific formulation we give is taken from \cite{alspach1977primariness}.

\begin{theorem}[Classification theorem of ordinal spaces]\label{th: classification-ordinals} Let $\alpha < \beta$ be two ordinals of the same cardinality, and let $\xi$ be the first ordinal of this cardinality. Write $\alpha = \xi\eta_1 + \rho_1$ and $\beta = \xi\eta_2 + \rho_2$ ($\rho_1 < \xi,\rho_2 < \xi$). Then there are two cases:
\begin{enumerate}
    \item \label{it: ord1} If the cardinality of $\alpha$ is an uncountable regular cardinal and $\eta_1 \leq \xi$, then $C([0, \alpha])$ is isomorphic to $C([0, \beta])$ iff $\eta_1$ and $\eta_2$ have the same cardinality.
    \item \label{it: ord2} Otherwise, $C([0, \alpha])$ is isomorphic to $C([0, \beta])$ iff $\beta < \alpha^\omega$ (where $\omega$ is the first infinite ordinal).
\end{enumerate}
\end{theorem}

\begin{proof}[Proof of Theorem \ref{th: Ogden}]
    We start with the observation that $C([0, \alpha2])$ is naturally isomorphic to $C([0, \alpha])\oplus C([0, \alpha])$. It follows from Johnson's automatic continuity result \cite{johnson1967continuity} that whenever $C([0, \alpha2]) \sim C([0, \alpha])$, then all homomorphisms from $\mathscr{B}(C([0, \alpha]))$ into a Banach algebra are continuous. By the Classification Theorem \ref{th: classification-ordinals}, the only cases when this does not hold (up to isomorphism of the space of continuous functions) is when $\alpha$ is finite or $\alpha = \xi n$ where $\xi$ is the first ordinal with the same cardinality as $\alpha$, the cardinality of $\alpha$ is an uncountable regular cardinal and $n \in \mathbb{N}$. The finite case is obvious, while in the latter case $C([0, \alpha]) \sim C([0, \xi])^n$ and Ogden's theorem \cite{ogden1996homomorphisms} gives that every homomorphism from $\mathscr{B}(C([0, \xi]))$ to a Banach algebra is continuous. The result now follows from Corollary \ref{cor: maximality norm}.
\end{proof}

\vspace{1em}

\noindent\textbf{Acknowledgements.} This paper is part of the author's PhD research at Lancaster University under the supervision of Professor N.~J.~Laustsen. The author is deeply indebted to Professor Laustsen for both introducing the problem and providing guidance throughout its resolution; his insightful comments and meticulous revisions have had a significant impact on this work. He also extends his gratitude to Dr.~Max Arnott and Dr.~Matthew Daws for their valuable comments and suggestions on earlier versions of the manuscript, and to Professor W.~B.~Johnson for kindly explaining how the complex version of Theorem \ref{th: hagler} follows and for permitting the inclusion of his argument here. Finally, the author thanks the anonymous referee for an exceptionally detailed report that significantly improved the clarity and overall presentation of the paper.

He also acknowledges with thanks the funding from the EPSRC (grant number EP/W524438/1) that has supported his studies.


\end{document}